\documentclass[10pt]{amsart}
\usepackage{amssymb}

\theoremstyle{definition}
\newtheorem{thm}{Theorem}[section]
\newtheorem{lem}[thm]{Lemma}
\newtheorem{prp}[thm]{Proposition}

\newtheorem{cor}[thm]{Corollary}

\newtheorem{rmk}[thm]{Remark}

\newtheorem{pbm}[thm]{Problem}

\newenvironment{pff}{{\em Proof:}}{\QED}

\newcommand{\beq}{\begin{equation}}
\newcommand{\eeq}{\end{equation}}
\newcommand{\beqr}{\begin{eqnarray*}}
\newcommand{\eeqr}{\end{eqnarray*}}

\newcommand{\bit}{\begin{itemize}}
\newcommand{\eit}{\end{itemize}}
\newcommand{\ts}[1]{{\textstyle{#1}}}
\newcommand{\ds}[1]{{\displaystyle{#1}}}
\newcommand{\ssum}[2]{{\ts{ {\ds{\sum}}_{#1}^{#2} }}}
\newcommand{\sbigcup}[2]{{\ts{ {\ds{\bigcup}}_{#1}^{#2} }}}

\newcommand{\ep}{\varepsilon}
\newcommand{\zt}{\zeta}
\newcommand{\et}{\eta}

\newcommand{\sm}{\sigma}

\newcommand{\ph}{\varphi}

\newcommand{\om}{\omega}
\newcommand{\ta}{\tau}

\newcommand{\R}{{\mathbf{R}}}
\newcommand{\C}{{\mathbf{C}}}
\newcommand{\N}{{\mathbf{N}}}

\newcommand{\op}{^{\mathrm{op}}}

\newcommand{\tsr}{{\mathrm{tsr}}}

\newcommand{\dist}{{\mathrm{dist}}}

\newcommand{\RR}{{\mathrm{RR}}}

\newcommand{\sa}{{\mathrm{sa}}}

\newcommand{\dirlim}{\displaystyle \lim_{\longrightarrow}}

\newcommand{\andeqn}{\,\,\,\,\,\, {\mbox{and}} \,\,\,\,\,\,}
\newcommand{\QED}{\rule{0.4em}{2ex}}

\newcommand{\ca}{C*-algebra}

\newcommand{\ct}{continuous}
\newcommand{\mvn}{Murray-von Neumann equivalent}

\newcommand{\pj}{projection}
\newcommand{\Wolog}{Without loss of generality}

\title[A simple C*-algebra not isomorphic to its opposite]{A simple
separable C*-algebra not isomorphic to its opposite algebra}

\author{N.\  Christopher Phillips}

\address{Department of Mathematics, University  of Oregon,
       Eugene OR 97403-1222, USA.}

\date{22 Feb.\  2002}

\email[]{ncp@darkwing.uoregon.edu}

\thanks{Research partially supported by NSF grant DMS 0070776.}

\begin{document}

\setcounter{section}{-1}

\begin{abstract}
We give an example of a simple
separable C*-algebra which is not isomorphic to its opposite algebra.
Our example is nonnuclear and stably finite,
has real rank zero and stable rank one,
and has a unique tracial state.
It has trivial $K_1$, and its $K_0$-group is order isomorphic to a
countable subgroup of $\R$.
\end{abstract}

\maketitle

\section{Introduction}

\indent
The purpose of this note is to give an example of a simple
separable C*-algebra which is not isomorphic to its opposite algebra.
By the opposite algebra $A\op$ of a \ca\  $A,$ we mean the algebra $A$
with the multiplication reversed but all other operations, including
the scalar multiplication, the same.
(The opposite algebra is isomorphic to the
complex conjugate algebra, via the map $x \mapsto x^*$.)
The existence of type~I C*-algebras not isomorphic to their opposites
has been known for some time; early examples are due to Raeburn
and P.\  Green, and several examples with additional interesting
properties are given in \cite{PhX}.
It has been known for some time that there are von Neumann factors,
with separable preduals,
of type~II$_1$~\cite{Co2} and type~III~\cite{Co1} which are not
isomorphic as von Neumann algebras to their opposites.
A \ca\  isomorphism of von Neumann algebras is necessarily a
von Neumann algebra isomorphism, by Corollary~5.13 of~\cite{SZ}, so
these are simple C*-algebras not isomorphic to their opposite algebras.
However, one wants separable examples.

We construct our example by applying a method of Blackadar~\cite{Bl0}
to the type II$_1$ factor of Corollary~7 of~\cite{Co2}.
The resulting \ca\  is nonnuclear and stably finite,
has real rank zero~\cite{BP} and stable rank one~\cite{Rf},
and has a unique tracial state.
It has trivial $K_1$, and its $K_0$-group is order isomorphic to a
countable subgroup of $\R$.
However, we have little control over other properties.
In particular, we can't specify which subgroups of $\R$ occur,
although we can show that there are uncountably many of them.

The recent work on classification of simple nuclear \ca s,
for example \cite{Ph} and \cite{Kr} in the purely infinite case
and \cite{EGL} and \cite{Ln} in the stably finite case,
suggests that all simple nuclear \ca s might
be isomorphic to their opposites.
The algebras $A$ and $A\op$ always have the same Elliott invariant,
and the nonisomorphism $A \not\cong A\op$ in our example
shows one way in which the
Elliott conjecture goes wrong when the nuclearity condition is dropped.
Other examples of nonisomorphic simple separable nonnuclear \ca s
with the same Elliott invariant are known.
The algebras have been distinguished by the Haagerup invariant
(\cite{PhY}; proof of Theorem~4.3.8 of~\cite{Ph}),
which finite dimensional operator spaces can be embedded in the algebra
(\cite{PhY}; proof of Theorem~4.3.11 of~\cite{Ph}),
quasidiagonality (\cite{GP}, \cite{PhY}),
approximate divisibility
(Theorem~1.4 of~\cite{DR}; also see Remark~4.3.2 of~\cite{Ph}),
and tensor indecomposability of an associated von Neumann algebra
(\cite{GP}).
None of these methods is capable of distinguishing a \ca\  from its
opposite algebra.

% The example is also a simple unital \ca\  which is not
% isomorphic to the complexification of any real \ca.

\section{Blackadar's result and some analogs}

\indent
A key ingredient of our construction is the following result of
Blackadar, Proposition~2.2 of~\cite{Bl0}.
% It will be applied with the large algebra taken to be the type~II$_1$
% factor of Corollary~7 of~\cite{Co2}.

\begin{lem}\label{Simple}
Let $N$ be a simple \ca,
and let $A \subset N$ be a separable C* subalgebra.
Then there exists a simple separable C* subalgebra $B$ with
$A \subset B \subset N$.
\end{lem}

To obtain the other properties necessary for our construction,
we need to know that it is possible to find separable intermediate
subalgebras preserving other properties from the large algebra.
To just prove the existence of a separable simple \ca\  not
isomorphic to its opposite, we only need the next lemma, on traces.
The remaining lemmas will be used to show that the algebra can be
chosen to have additional good properties.
Some are already implicit in previous work.

\begin{lem}\label{Tr}
Let $N$ be a unital \ca,
and let $A \subset N$ be a separable C* subalgebra.
Then there exists a separable C* subalgebra $B$ with
$A \subset B \subset N$ such that every tracial state on $B$ is the
restriction of a tracial state on $N$.
\end{lem}

\begin{pff}
For any \ca\  $D$, let $[D, D]$ denote the linear
span of the commutators $[a, b] = a b - b a$ with $a, \, b \in D$.
Also, we let, for any $d \in D$ and $S \subset D$,
\[
\dist (d, S) = \inf \{ \| d - x \| \colon x \in S \}.
\]

\Wolog\  $A$ contains the identity of $N$.
We construct inductively separable C* subalgebras
$B_n \subset N$ such that:
\begin{itemize}
\item
$B_0 = A$.
\item
$B_0 \subset B_1 \subset B_2 \subset \cdots$.
\item
For every $b \in B_n$, we have
$\dist (b, \, [B_{n + 1}, \, B_{n + 1} ] ) = \dist (b, \, [N, N] )$.
\end{itemize}
We do the induction step; the base case is the same.
Given $B_n$, choose a countable dense subset $S \subset B_n$.
For $b \in S$ and $m \in \N$, choose
\[
l (b, m) \in \N \andeqn
y_{b, m, 1}, \, y_{b, m, 2}, \, \dots, \, y_{b, m, l (b, m)}, \,
z_{b, m, 1}, \, z_{b, m, 2}, \, \dots, \, z_{b, m, l (b, m)} \in N
% y_{b, m, 1}, \, \dots, \, y_{b, m, l (b, m)}, \,
% z_{b, m, 1}, \, \dots, \, z_{b, m, l (b, m)} \in N
\]
such that
\[
\left\| b -
 \ssum{j = 1}{l (b, m)} [y_{b, m, j}, \, z_{b, m, j}] \right\|
< \frac{1}{m} + \dist (b, \, [N, N] ).
\]
Take $B_{n + 1}$ to be the separable C* subalgebra of $N$
generated by $B_n$ and the countable set
\[
\{ y_{b, m, j}, \, z_{b, m, j} \colon
  {\mbox{$b \in S$, $m \in \N$, and $1 \leq j \leq l (b, m)$}} \}.
\]
For $c \in B_n$ and $\ep > 0$, choose $b \in S$ with
$\| c - b \| < \frac{1}{3} \ep$, choose $m \in \N$ with
$\frac{1}{m} < \frac{1}{3} \ep$, and set
\[
s = \sum_{j = 1}^{l (b, m)} [y_{b, m, j}, \, z_{b, m, j}]
  \in [B_{n + 1}, \, B_{n + 1} ].
\]
Then
$\dist (c, \, [N, N] )
  < {\textstyle{ \frac{1}{3} }} \ep + \dist (b, \, [N, N] )$,
whence
\[
\| c - s \| \leq \| c - b \| + \| b - s \|
 < {\textstyle{ \frac{1}{3} }} \ep
  + \left[ {\textstyle{ \frac{1}{m} }} + \dist (b, \, [N, N] ) \right]
 < \ep + \dist (c, \, [N, N] ).
\]
% \begin{align*}
% \| c - s \| & \leq \| c - b \| + \| b - s \|
%  < {\textstyle{ \frac{1}{3} }} \ep
%   + \left[ {\textstyle{ \frac{1}{m} }} + \dist (b, \, [N, N] ) \right]
%                    \\
% & < \ep + \dist (c, \, [N, N] ).
% \end{align*}
Since $\ep > 0$ is arbitrary, this gives
\[
\dist (c, \, [B_{n + 1}, \, B_{n + 1} ] ) \leq \dist (c, \, [N, N] ),
\]
which completes the induction step.

Now set
\[
B = {\overline{\sbigcup{n = 0}{\infty}  B_n}}.
\]
It is clear that for every $b \in \bigcup_{n = 0}^{\infty} B_n$,
\[
\dist (b, \, [B, B] ) = \dist (b, \, [N, N] ),
\]
and equality easily follows for all $b \in B$.
It is now immediate that
\[
\dist (b, \, {\overline{[B, B]}} ) = \dist (b, \, {\overline{[N, N]}} )
\]
for all $b \in B$.
This implies that the inclusion of $B$ in $N$ defines an isometric
linear map
\[
T \colon B / {\overline{[B, B]}} \to N / {\overline{[N, N]}}.
\]

Let $\ta \colon B \to \C$ be any tracial state.
We construct a tracial state $\sm$ on $N$ such that $\sm |_B = \ta$.
By continuity and the trace property, $\ta$ induces a linear
functional ${\overline{\ta}} \colon B / {\overline{[B, B]}} \to \C$
with $\| {\overline{\ta}} \| = 1$.
The Hahn-Banach Theorem provides a linear functional
$\om \colon N / {\overline{[N, N]}} \to \C$ such that
$\om \circ T = {\overline{\ta}}$ and $\| \om \| = 1$.
Let $\sm \colon N \to \C$ be the composition of $\om$ with the
quotient map $N \to N / {\overline{[N, N]}}$.
Then $\sm |_B = \ta$, and in particular $\sm (1) = 1$.
Since $\| \sm \| = 1$, it follows that $\sm$ is a state.
Moreover, $\sm$ is a trace because it vanishes on $[N, N]$.
So $\sm$ is the required tracial state.
\end{pff}

\begin{lem}\label{StRank}
Let $N$ be a unital \ca,
and let $A \subset N$ be a separable C* subalgebra.
Then there exists a separable C* subalgebra $B$ with
$A \subset B \subset N$ such that $\tsr (A) \leq \tsr (N)$.
\end{lem}

\begin{pff}
\Wolog\  $A$ contains the identity of $N$.
Following Definition~1.4 and Proposition~1.6 of~\cite{Rf},
we let $r = \tsr (N)$ and we construct $B$ in such a way that the space
${\mathrm{Lg}}_r (B) \subset B^r$ (Notation~1.3 of~\cite{Rf}),
consisting of all $b = (b_1, b_2, \dots, b_r) \in B^n$ such that
$\{b_1, b_2, \dots, b_r\}$ generates $B$ as a left ideal,
is dense in $B^r$.

We construct inductively separable C* subalgebras
$B_n \subset N$ such that:
\begin{itemize}
\item
$B_0 = A$.
\item
$B_0 \subset B_1 \subset B_2 \subset \cdots$.
\item
$B_n^r \subset {\overline{ {\mathrm{Lg}}_r (B_{n + 1} ) }}$ for all $n$.
\end{itemize}
We do the induction step; the base case is the same.
Choose a suitable norm on $N^r$.
Let $S$ be a countable dense subset of $B_n^r$.
For each $b \in S$ and $m \in \N$, use $\tsr (N) = r$ to choose
\[
x_{b, m} = (x_{b, m, 1}, \, x_{b, m, 2}, \, \dots, \, x_{b, m, r})
  \in {\mathrm{Lg}}_r (N)
\]
such that $\| x_{b, m} - b \| < 1/m$.
By definition, there are
\[
y_{b, m, 1}, \, y_{b, m, 2}, \, \dots, \, y_{b, m, r} \in N
\]
such that
\[
y_{b, m, 1} x_{b, m, 1} + y_{b, m, 2} x_{b, m, 2}
  + \cdots + y_{b, m, r} x_{b, m, r} = 1.
\]
Take $B_{n + 1}$ to be the separable C* subalgebra of $N$
generated by $B_n$ and the countable set
\[
\{ x_{b, m, j}, \, y_{b, m, j} \colon
  {\mbox{$b \in S$, $m \in \N$, and $1 \leq j \leq r$}} \}.
\]
Clearly each $x_{b, m}$ is in ${\mathrm{Lg}}_r (B_{n + 1})$,
so the induction step is complete.

Now set
\[
B = {\overline{\sbigcup{n = 0}{\infty}  B_n}}.
\]
We have
\[
B^r = {\overline{\sbigcup{n = 1}{\infty}  B_{n - 1}^r }}
  \subset {\overline{\sbigcup{n = 1}{\infty} {\mathrm{Lg}}_r (B_n) }}
  \subset {\overline{ {\mathrm{Lg}}_r (B) }},
\]
whence $\tsr (B) \leq r$.
\end{pff}

\begin{lem}\label{RRank}
Let $N$ be a \ca, and let $A \subset N$ be a separable C* subalgebra.
Then there exists a separable C* subalgebra $B$ with
$A \subset B \subset N$ such that $\RR (A) \leq \RR (N)$.
\end{lem}

\begin{pff}
The proof is the same as for Lemma~\ref{StRank}.
Following \cite{BP}, we simply adjust the indexing and consider
only ${\mathrm{Lg}}_{r + 1} (B) \cap (B_{\sa})^{r + 1}$, etc.
\end{pff}

\begin{lem}\label{InjOnKThy}
Let $N$ be a \ca, and let $A \subset N$ be a separable C* subalgebra.
Then there exists a separable C* subalgebra $B$ with
$A \subset B \subset N$ such that the map $K_* (B) \to K_* (N)$
is injective and induces an order isomorphism of $K_0 (B)$ with
a subgroup of $K_0 (N)$.
\end{lem}

\begin{pff}
Unitizing, we may assume that $N$ is unital
and $A$ contains the identity of $N$.

If $D$ is a unital \ca,
we write $u \sim_D v$ for unitaries $u, \, v$ in some matrix algebra
$M_m (D)$ which are homotopic in the unitary group $U (M_m (D))$.
We also write $p \sim_D q$ for \pj s $p, \, q \in M_m (D)$ which are
\mvn\  in $M_m (D)$, and $p \precsim_D q$ if $p$ is \mvn\  to a
sub\pj\  of $q$ in $M_m (D)$.
Finally, we write $1_m$ for the identity of $M_m (D)$.

We construct inductively separable C* subalgebras
$B_n \subset N$ such that:
\begin{itemize}
\item
$B_0 = A$.
\item
$B_0 \subset B_1 \subset B_2 \subset \cdots$.
\item
For $m \in \N$ and $u \in U (M_m (B_n))$,
if $u \sim_N 1_m$ then $u \sim_{B_{n + 1}} 1_m$.
\item
For $m \in \N$ and \pj s $p, \, q \in M_m (B_n)$,
if $p \sim_N q$ then $p \sim_{B_{n + 1}} q$.
\item
For $m \in \N$ and \pj s $p, \, q \in M_m (B_n)$,
if $p \precsim_N q$ then $p \precsim_{B_{n + 1}} q$.
\end{itemize}

We do the induction step; the base case is the same.
Thus, suppose that $B_n$ has been found.
Since two unitaries $u$ and $v$ with $\| u - v \| < 2$ are homotopic,
and since $B_n$ is separable, for each $m$ there are only countably
many homotopy classes $[u]$ of unitaries $u \in M_m (B_n)$.
Let $S_0$ be the set of all pairs $(m, [u])$ with $m \in \N$
and $u \in U (M_m (B_n))$ such that $u \sim_N 1_m$.
Then $S_0$ is countable.
For each $(m, [u]) \in S_0$, choose a unitary path $t \mapsto w (t)$
in $M_m (N)$ with $w_0 = u$ and $w_1 = 1$, choose
$0 = t_0 < t_1 < \cdots < t_k = 1$ such that
$\| w (t_j) - w (t_{j - 1}) \| < 2$ for $1 \leq j \leq k$, and let
$T_0^{ (m, [u]) }$ be the subset of $N$ consisting of all
matrix entries of all $w (t_j)$.
Then let $T_0$ be the countable set
\[
T_0 = \bigcup_{(m, [u]) \in S_0} T_0^{ (m, [u]) }.
\]
It is easy to see from the transitivity of homotopy
that whenever $u \in U (M_m (B_n))$ satisfies $u \sim_N 1_m$,
then also $u \sim_D 1_m$ for any \ca\  $D$ containing $B_n$ and $T_0$.

Since \pj s $p$ and $q$ with $\| p - q \| < 1$ are \mvn,
each $M_m (B_n)$ contains only countably many
Murray-von Neumann equivalence classes of \pj s.
A similar construction produces a countable subset $T_1 \subset N$
such that whenever $p, \, q \in M_m (B_n)$ are \pj s which are
\mvn\  in $M_r (N)$ but not in $M_m (B_n)$, then there exist
\pj s $p_0, \, q_0 \in M_m (B_n)$ and a matrix $s \in M_m (N)$
such that $p_0 \sim_{B_n} p$ and $q_0 \sim_{B_n} q$,
such that all the entries of $s$ are in $T_1$,
and such that $s^* s = p_0$ and $s s^* = q_0$.
It follows that whenever \pj s $p, \, q \in M_m (B_n)$ satisfy
$p \sim_N q$, then $p \sim_D q$ for any
\ca\  $D$ containing $B_n$ and $T_1$.
By essentially the same method, one can construct a countable
subset $T_2 \subset N$ such that whenever $D$ is a
\ca\  containing $B_n$ and $T_0$, and whenever
$p, \, q \in M_m (B_n)$ are \pj s such that
$p \precsim_N q$, then $p \precsim_D q$.
The induction step is now completed by taking $B_{n + 1}$ to be the
C* subalgebra of $N$ generated by $B_n \cup T_0 \cup T_1 \cup T_2$.

Now set
\[
B = {\overline{\sbigcup{n = 0}{\infty}  B_n}}.
\]
For every $m$ and every $u \in U (M_m (B))$ such that $u \sim_N 1_m$,
there is $n$ and $v \in U (M_m (B_n))$ such that $v \sim_B u$.
So $v \sim_{B_{n + 1}} 1_m$, whence $u \sim_B 1_m$.
This implies that $K_1 (B) \to K_1 (N)$ is injective.
By a similar argument, for every $m$ and for any two \pj s
$p, \, q \in M_m (B_n)$ such that $p \sim_N q$, we have $p \sim_B q$.
Therefore $K_0 (B) \to K_0 (N)$ is injective.
It remains to prove that $K_0 (B) \to K_0 (N)$ is
an order isomorphism onto its image.
Since this map preserves order, we need only show that if
$\et \in K_0 (B)$ is a class whose image in $K_0 (N)$ is positive,
then $\et > 0$.
So let $p, \, q$ be \pj s in matrix algebras over $B$ such that
$\et = [q] - [p]$, and let $e$ be a \pj\  in some matrix algebra over
$N$ such that $[q] - [p] = [e]$ in $K_0 (N)$.
\Wolog\  there is $n$ such that $p$ and $q$
are in matrix algebras over $B_n$.
Replacing $p$ and $q$ by $p \oplus 1_r$ and $q \oplus 1_r$
(which are still in matrix algebras over $B_n$) for
suitable $r$, we may assume that in a suitable matrix algebra
$M_m (N)$, we can write $q = p_0 + e_0$ with
$p_0 \sim_N p$ and $e_0 \sim_N e$.
In particular, $p \precsim_N q$.
By construction $p \precsim_{B_{n + 1}} q$, whence $p \precsim_B q$.
It follows that $\et = [q] - [p] > 0$ in $K_0 (B)$.
\end{pff}

\section{The main result}

\indent
Our main result will follow from the following proposition,
using a suitable choice of the type~II$_1$ factor.

\begin{prp}\label{Combination}
Let $N$ be a type~II$_1$ factor with separable predual and
with trace $\ta$.
Let $G_0$ be a countable subgroup of $K_0 (N)$, which we identify
with $\R$ via $\ta_*$.
Then there exists a simple separable unital
weak operator dense C* subalgebra $A \subset N$
such that $\tsr (A) = 1$, $\RR (A) = 0$, $K_1 (A) = 0$,
the map $K_0 (A) \to K_0 (N)$ induces an order isomorphism of
$K_0 (A)$ with a subgroup of $K_0 (N)$ containing $G_0$,
and $A$ has as unique tracial state the restriction $\ta |_A$.
\end{prp}

\begin{pff}
Let $S \subset N$
be a countable subset which is weak operator dense in $N$.
For each $g \in G_0$ with $g > 0$, choose an integer $m (g) > 0$ and
a \pj\  $p_g \in M_{m (g)} (N)$ such that $\ta (p_g) = g$.
Let $P \subset N$ be the unital C* subalgebra of $N$
generated by $S$ and all
the matrix entries of all $p_g$ for $g \in G_0 \cap (0, \infty)$.
Then $P$ is separable and the image of the map
$K_0 (P) \to K_0 (N)$ contains $G_0$.

We now construct by induction on $n$ separable subalgebras
$A_n$, $B_n$, $C_n$, $D_n$, and $E_n$ with
\[
P \subset A_0 \subset B_0 \subset C_0
                          \subset D_0 \subset E_0 \subset \cdots
  \subset A_n \subset B_n \subset C_n
                          \subset D_n \subset E_n \subset \cdots
\]
and such that, for all $n$, we have:
\begin{itemize}
\item
$A_n$ is simple.
\item
$B_n$ has as unique tracial state the restriction $\ta |_{A_0}$.
\item
$\tsr (C_n) = 1$.
\item
$\RR (D_n) = 0$.
\item
The map $K_0 (E_n) \to K_0 (N)$ is an order isomorphism onto its
image and the map $K_1 (E_n) \to K_1 (N)$ is injective.
\end{itemize}

The base case is the same as the induction step, using
$P$ in place of $E_n$, so we do only the induction step.
Suppose the subalgebras have been constructed through $E_n$.
Use Lemma~\ref{Simple} to choose a simple separable C* subalgebra
$A_{n + 1}$ with $E_n \subset A_{n + 1} \subset N$.
Use Lemma~\ref{Tr} to choose a separable C* subalgebra
$B_{n + 1}$ with $A_{n + 1} \subset B_{n + 1} \subset N$
such that every tracial state on $B_{n + 1}$ is the
restriction of a tracial state on $N$.
The invertible elements in $N$ are dense, since in the polar
decomposition $x = s (a^* a)^{1/2}$ of any $x \in N$ we can first
replace the partial isometry by a unitary and then, with an error of
$\ep$, replace $(a^* a)^{1/2}$ by $(a^* a)^{1/2} + \ep \cdot 1$.
So we can use Lemma~\ref{StRank} to choose a
separable C* subalgebra $C_{n + 1}$ with
$B_{n + 1} \subset C_{n + 1} \subset N$
such that $\tsr (C_{n + 1}) = 1$.
Use Lemma~\ref{RRank} and $\RR (N) = 0$ to choose a separable
C* subalgebra $D_{n + 1}$ with $C_{n + 1} \subset D_{n + 1} \subset N$
such that $\RR (D_{n + 1}) = 0$.
Use Lemma~\ref{InjOnKThy} to choose a separable C* subalgebra
$E_{n + 1}$ with $D_{n + 1} \subset E_{n + 1} \subset N$
such that the map $K_1 (E_{n + 1}) \to K_1 (N)$ is injective
and such that the map $K_0 (E_{n + 1}) \to K_0 (N)$ induces an
order isomorphism onto its image.

Now set
\[
A = {\overline{\sbigcup{n = 0}{\infty} A_n }}
  = {\overline{\sbigcup{n = 0}{\infty} B_n }}
  = {\overline{\sbigcup{n = 0}{\infty} C_n }}
  = {\overline{\sbigcup{n = 0}{\infty} D_n }}
  = {\overline{\sbigcup{n = 0}{\infty} E_n }}.
\]
We verify that $A$ has the required properties.
Obviously $A$ is separable.
The algebra $A$ is weak operator dense in $N$ because it contains $S$.
{}From $A = {\overline{\bigcup_{n = 0}^{\infty} A_n }}$ and simplicity
of the $A_n$, a standard argument shows that $A$ is simple.
Any trace $\ta_0$ on $A$ must restrict to a trace on each $B_n$,
necessarily $\ta |_{B_n}$.
Since $\bigcup_{n = 0}^{\infty} B_n$ is dense in $A$, it follows
that $\ta_0 = \ta |_A$.
On the other hand, clearly $\ta |_A$ is a trace on $A$.
We have $\tsr (A) = 1$ by Theorem~5.1 of~\cite{Rf}, because
$A = \dirlim C_n$ and $\tsr (C_n) = 1$ for all $n$.
It is clear that $\RR (A) = 0$, because
$A = \dirlim D_n$ and $\RR (D_n) = 0$ for all $n$.
Finally, using the relation $A = \dirlim E_n$, a slightly easier
version of the argument of the
last paragraph of the proof of Lemma~\ref{InjOnKThy} shows that
the map $K_1 (A) \to K_1 (N)$ is injective and
the map $K_0 (A) \to K_0 (N)$ is an order isomorphism onto its image.
Since $K_1 (N) = 0$, this immediately gives $K_1 (A) = 0$.
Moreover, since $P \subset A$ it follows that $G_0$ is contained in
the image of $K_0 (A)$.
\end{pff}

\begin{thm}\label{Main}
Let $G_0$ be a countable subgroup of $\R$.
Then there exists a simple separable stably finite
unital C* subalgebra $A$ with $A \not\cong A\op$ and
such that $A$ has stable rank one and real rank zero, $K_1 (A) = 0$,
the group $K_0 (A)$ is isomorphic as a scaled ordered group to a
countable subgroup of $\R$ containing $G_0$,
and $A$ has a unique tracial state.
\end{thm}

\begin{pff}
Let $N$ be the type II$_1$ factor of Corollary~7 of~\cite{Co2},
which is not isomorphic as a von Neumann algebra to $N\op$.
Apply Proposition~\ref{Combination} with this $N$ and with $G_0$ as in
the hypotheses, and let $A$ be the resulting \ca.
The only property that is not immediate is the
nonisomorphism $A \not\cong A\op$.

Suppose that there is an isomorphism $\ph \colon A \to A\op$.
Let $\ta_0$ be the unique trace on $A$, and let $\ta_0\op$ be $\ta_0$
regarded as a trace on $A\op$.
Let $\pi$ and $\pi\op$ be the Gelfand-Naimark-Segal representations
of $A$ and $A\op$ associated with $\ta_0$ and $\ta_0\op$.
Then $\ta_0 = \ta_0\op \circ \ph$ by uniqueness of the traces, whence
$\pi$ is unitarily equivalent to $\pi\op \circ \ph$.
It follows that $\pi\op (A\op)''$ is isomorphic as a von Neumann
algebra to $\pi (A)''$.

We claim that $\pi (A)'' \cong N$.
Let $\ta$ be the trace on $N$, and note that
Proposition~\ref{Combination} gives $\ta |_A = \ta_0$.
We may assume that $N$ is represented in the canonical way on
the Hilbert space $L^2 (N, \ta)$; this is just the
Gelfand-Naimark-Segal representation of $N$ associated with $\ta$.
All we need to know about it is contained in Proposition~III.3.12
of~\cite{Tk}.
The Hilbert space $H_{\pi}$ of the representation $\pi$ is
by construction a subspace of $L^2 (N, \ta)$.
We show below that the \ca\  $A$ is dense in $N$ in the norm
$\| a \|_2 = \ta (a^* a)^{1/2}$ associated with $L^2 (N, \ta)$.
Since $A$ is contained in $H_{\pi}$,
it follows that $H_{\pi} = L^2 (N, \ta)$.
Therefore, using weak operator density of $A$ in $N$,
we get $\pi (A)'' = N$, proving the claim.

To prove the density statement, let $x \in N$.
Write $x = a + i b$ with $a, \, b \in N_{\sa}$.
Using the Kaplansky Density Theorem, Theorem~2.3.3 of~\cite{Pd1},
and the fact that the strong operator topology is metrizable on
bounded sets when the predual is separable, find bounded sequences
$(a_n)$ and $(b_n)$ in $A_{\sa}$ such that $a_n \to a$ and $b_n \to b$
in the strong operator topology.
Set $x_n = a_n + i b_n \in A$.
Then $x_n \to x$ and $x_n^* \to x^*$ in the strong operator topology.
Since multiplication is jointly strong operator \ct\  on bounded sets,
it follows that $(x_n - x)^* (x_n - x) \to 0$
in the strong operator topology.
The trace $\ta$ is strong operator \ct, so
$\| x_n - x \|_2^2 = \ta ( (x_n - x)^* (x_n - x) ) \to 0$.
This proves density.

Similarly $\pi\op (A\op)'' \cong N\op$.
But now we have contradicted the property $N \not\cong N\op$.
\end{pff}

\section{Consequences and open problems}

\indent
Recall that any real \ca\  $B$ has a complexification
$B_{\C} = B \otimes_{\R} \C$,
which is a complex \ca\  with a conjugate linear automorphism
$a \otimes \zt \mapsto a \otimes {\overline{\zt}}$.
(We refer to Part~II of \cite{Gd} for the general theory of
real \ca s.)
In particular, $B_{\C}\op \cong B_{\C}$.
Therefore we obtain the following corollary.

\begin{cor}\label{NotCx}
There exists a simple separable stably finite unital C* subalgebra $A$
which is not isomorphic to the complexification of any real \ca.
\end{cor}

\begin{rmk}\label{NonNuc}
The algebras in Theorem~\ref{Main} are not nuclear,
because they are weak operator dense in a factor $N$ which is not
hyperfinite.
One can certainly force them to be nonexact, and it seems unlikely
that the construction can be made to produce exact \ca s.
\end{rmk}

\begin{rmk}\label{Many}
Theorem~\ref{Main} produces uncountably many mutually nonisomorphic
examples, because
no countable union of countable subgroups of $\R$ can contain
all countable subgroups of $\R$.
%
% If there were only countably many, then the subgroup $G$ of $\R$
% generated by the union of the $K_0$-groups of these examples
% would be countable.
% A contradiction is obtained by taking $G_0$ to contain any
% element of $\R \setminus G$.
\end{rmk}

We close by giving several open problems.

\begin{pbm}\label{GenKThy}
Let $F$ be a simple unital AF algebra.
Find a simple separable stably finite unital C* subalgebra $A$
with $A \not\cong A\op$ and such that $\tsr (A) = 1$, $\RR (A) = 0$,
and $A$ has the same Elliott invariant as $F$.
\end{pbm}

As far as we can tell, the methods of \cite{PhY} will not work,
because there is no reason to think the property $A \not\cong A\op$
is preserved through the steps of the construction there.

\begin{pbm}\label{Nat}
Find a more natural example of a simple separable \ca\  $A$ with
$A \not\cong A\op$.
\end{pbm}

The obvious approach is to try a \ca ic version of the constructions
of \cite{Co2} or \cite{Co1}, both of which involve crossed products.

\begin{pbm}\label{PI}
Is there a purely infinite simple separable \ca\  $A$
such that $A \not\cong A\op$?
\end{pbm}

The methods here don't seem to apply to the infinite case.
One might hope that a \ca ic version of the construction
of \cite{Co1} could produce such an example.

\begin{pbm}\label{Exact}
Is there a separable exact \ca\  $A$ such that $A \not\cong A\op$?
\end{pbm}

\end{document}